\documentclass[11pt]{article}
\usepackage{amsmath,amssymb,amsbsy,amsfonts,latexsym,amsthm,mathrsfs}
\usepackage{graphicx}
\usepackage{float}
\usepackage{color}
\usepackage{todonotes}
\newtheorem{theorem}{Theorem}[section]
\newtheorem{proposition}[theorem]{Proposition}

\newtheorem{condition}[theorem]{Condition}

\newtheorem{lemma}[theorem]{Lemma}
\newtheorem{remark}[theorem]{Remark}
\newtheorem{example}[theorem]{Example}
\begin{document}
\title{Asymptotic results for compound sums in separable Banach spaces}
\author{Claudio Macci\thanks{Address: Dipartimento di Matematica,
Universit\`a di Roma Tor Vergata, Via della Ricerca Scientifica,
I-00133 Rome, Italy. e-mail: \texttt{macci@mat.uniroma2.it}} \and
Barbara Pacchiarotti\thanks{Address: Dipartimento di Matematica,
Universit\`a di Roma Tor Vergata, Via della Ricerca Scientifica,
I-00133 Rome, Italy. e-mail: \texttt{pacchiar@mat.uniroma2.it}}}
\maketitle
\begin{abstract}
	\noindent We prove large and moderate deviation results for sequences of compound sums,
	where the summands are i.i.d. random variables taking values in a separable Banach space.
	We establish that the results hold by proving that we are dealing with exponentially tight sequences. We 
	present two moderate deviation results: in the first one the summands are centered, in the 
	second one the compound sums are centered.\\
	
\end{abstract}

\section{Introduction}
Large deviations give an asymptotic computation of small probabilities on an exponential scale (see, e.g., \cite{DZ:97} as a reference 
on this topic). We recall some basic definitions (see, e.g., Section 1.2 in \cite{DZ:97}; see also the pioneering paper
of \cite{Varadhan}). Throughout this paper a speed function is a sequence $\{v_n:n\geq 1\}$ such that $\lim_{n\to\infty}v_n=\infty$.
Moreover, given a sequence of random variables $\{Z_n:n\geq 1\}$ defined on the same probability space 
$(\Omega,\mathcal{F},P)$ and taking values in a topological space $\mathcal{X}$, the sequence $\{Z_n:n\geq 1\}$ satisfies
the large deviation principle (LDP for short) with rate function $I$ and speed function $v_n$ if: the function $I:\mathcal{X}\to[0,\infty]$
is a lower semicontinuous,
$$\liminf_{n\to\infty}\frac{1}{v_n}\log P(Z_n\in O)\geq-\inf_{x\in O}I(x)\quad\quad\mbox{for all open sets $O$,}$$
and
$$\limsup_{n\to\infty}\frac{1}{v_n}\log P(Z_n\in C)\leq-\inf_{x\in C}I(x)\quad\quad\mbox{for all closed sets $C$.}$$
Obviously the rate function $I$ cannot be equal to infinity everywhere.

The concept of exponential tightness plays a crucial role in large deviations; in fact this condition is often required to establish
that the large deviation principle holds for a sequence of random variables taking values in an infinite dimensional topological 
vector space. For completeness we recall that the sequence of random variables $\{Z_n:n\geq 1\}$ is exponentially tight with respect 
to the speed function $v_n$ if, for all $b>0$, there exists a compact $K_b\subset\mathcal{X}$ such that
$$\limsup_{n\to\infty}\frac{1}{v_n}\log P(Z_n\notin K_b)\leq-b.$$
Then we can say that, if we have the weak LDP (i.e. the lower bound for open sets, and the upper bound for the compact
sets only) and the exponential tightness, the full LDP holds with a good rate function (see, e.g., Lemma 1.2.18 in \cite{DZ:97}). Moreover,
in some situations (for instance locally compact spaces or Polish spaces), the exponential tightness is implied by the full LDP with a good 
rate function (see, e.g., Remark (a) at page 8 in \cite{DZ:97}).

In this paper we consider sequences of random variables taking values in an infinite dimensional topological vector space. 
Indeed we consider sequences of compound sums $\{S_{N_n}:n\geq 1\}$, where
$$S_n:=X_1+\cdots+X_n\quad\quad\mbox{for every}\ n\geq 1,$$
$\{X_n:n\geq 1\}$ are superexponential i.i.d. random variables taking values in a separable Banach space $B$ (see Condition \ref{cond:superexponential-tails-summands}), and independent of some integer-valued nonnegative random variables 
$\{N_n:n\geq 1\}$ which satisfy some suitable properties (see Condition \ref{cond:LD-*}).

For some technical reasons throughout this paper we assume that $B$ is separable. Alternatively one could assume 
some slightly weaker conditions: for instance the random variables $\{X_n:n\geq 1\}$ are separably valued, or the (common) 
distribution of those random variables is tight (and the random variables are defined on a complete probability space).

Compound sums have several applications in the literature; for instance, in insurance, the aggregate claim 
amount (up to a certain time $t$) is often modelled as a random sum of real valued random variables (see, e.g., \cite{SP_Ins_Fin}). Large 
and moderate deviations for compound sums have been already studied, essentially when $B=\mathbb{R}^h$ for some $h$ (see, e.g., Exercise 
2.3.19 in \cite{DZ:97} with some conditions on the common distribution of the random variables $\{X_n:n\geq 1\}$); other references for the 
case $h=1$ are \cite{EicheLo} and \cite{MitaTJM1997}. The moderate deviation results in \cite{EicheLo} are obtained by applying the 
Lindeberg's method invented to prove some central limit theorems. In particular some results in \cite{EicheLo} are derived by considering 
some bounds on cumulants (as in \cite{DoEiche} which concerns a quite general class of families of random variables). In \cite{EicheLo} the 
authors consider some appropriately centered and scaled compound sums that may converge to some Gaussian limiting distributions or some other 
non-Gaussian limiting distributions in (see, e.g., \cite{GneKo}, Chapters 3 and 4, for the non-Gaussian case). We also remark that Assumption 
4.2 in \cite{EicheLo} yields Conditions \ref{cond:superexponential-tails-summands} and \ref{cond:LD-*} in this paper with $B=\mathbb{R}$.

Our main motivation is to give a first contribution for compound sums in which the summands take values in an infinite 
dimensional space. In particular we are motivated by potential applications with random superpositions of random signals represented by 
continuous functions on a bounded interval; so, in this paper, we consider compound sums in which the summands take values in a separable 
Banach space $B$. So, for completeness, we also recall \cite{prochno2024large} as a recent reference (on large and moderate deviations) for
certain geometric random quantities in Banach spaces.

In this paper we prove a \emph{reference LDP} with speed $v_n=n$ (Proposition \ref{prop:LD}). Moreover we prove two results on 
\emph{moderate deviations}: the first one concerns the case of ``centered summands'' (Proposition \ref{prop:MD-summands-are-centered}), 
while the second one (Proposition \ref{prop:MD-centered-compound-sum}) concerns the case of ``centered compound sums''. 
In our moderate deviation results we consider the family of positive sequences (called \emph{scalings}) $\{a_n:n\geq 1\}$ such that
\begin{equation}\label{eq:MD-conditions}
	a_n\to 0\quad\mathrm{and}\quad na_n\to\infty;
\end{equation}
moreover some further hypotheses (see Condition \ref{cond:MD-*}) are needed.

We remark that, when we prove the exponential tightness conditions for large and moderate deviations, we adapt the proofs of the same 
results for non-compound sums in the literature (see \cite{deAcosta1} and \cite{deAcosta3}, respectively).

We conclude with the outline of the paper. We start with some preliminaries in Section \ref{sec:preliminaries}. Large and moderate 
deviations results are presented in Sections \ref{sec:LD} and \ref{sec:MD}, respectively. Finally Section \ref{sec:examples} concerns
some examples.

\section{Preliminaries}\label{sec:preliminaries}
We start with the following condition.

\begin{condition}\label{cond:superexponential-tails-summands}
Let $B$ be a separable Banach space, with norm $\|\cdot\|$. Moreover, let $\{X_n:n\geq 1\}$ be i.i.d. $B$-valued 
random variables that satisfy the \emph{strong Cramér condition}, i.e. $\mathbb{E}[e^{t\|X_1\|}]<\infty$ for all $t>0$. 
Then let $\{S_n:n\geq 1\}$ be defined by $S_n:=X_1+\cdots+X_n$ for all $n\geq 1$. We recall that, since $B$ is 
separable, the (common) distribution of the random variables $\{X_n:n\geq 1\}$ is tight; see, e.g., Theorem 1.4 in \cite{Billingsley-Convergence-Prob-Measures-1968}.
\end{condition}

Moreover, as far as the separable Banach space $B$ in Condition \ref{cond:superexponential-tails-summands} are concerned,
we denote the dual space of $B$ by $B^*$, the null vector in $B$ by $\widehat{0}$, the null vector in $B^*$ by $\widehat{0}^*$, and 
the \emph{duality between $B^*$ and $B$} by $\langle\cdot,\cdot\rangle$. Then, if we refer to the random variables $\{X_n:n\geq 1\}$
in Condition \ref{cond:superexponential-tails-summands}, we consider the function $\Lambda_X$ defined by
$$\Lambda_X(\theta):=\log\mathbb{E}[e^{\langle\theta,X_1\rangle}]\quad\quad\mbox{for all}\ \theta\in B^*,$$
and the Legendre transform of $\Lambda_X$, i.e. the function $\Lambda_X^*$ defined by
$$\Lambda_X^*(x):=\sup_{\theta\in B^*}\{\langle\theta,x\rangle-\Lambda_X(\theta)\}\quad\quad\mbox{for all}\ x\in B.$$
In particular (see, e.g., \cite{Ba-McK}), since we assume that $B$ is a separable Banach space (see Condition 
\ref{cond:superexponential-tails-summands}), we can say that:
\begin{equation}\label{eq:condition-for-the-mean}
	\mbox{there exists $\mu_X\in B$ such that $\langle v,\mu_X\rangle=\mathbb{E}[\langle v,X_1\rangle]$}\quad\quad (\mbox{for}\ v\in B^*);
\end{equation}
\begin{equation}\label{eq:condition-for-the-covariance}
	\left.\begin{array}{l}
		\mbox{there exists a linear operator $\Sigma_X:B^*\to B$ such that}\\
		\mbox{$\langle u,\Sigma_Xv\rangle=\langle\Sigma_Xu,v\rangle=\mathrm{Cov}(\langle u,X_1\rangle,\langle v,X_1\rangle)$}\quad\quad (\mbox{for}\ u,v\in B^*).
	\end{array}\right.
\end{equation}
We recall that $\mu_X$ in eq. \eqref{eq:condition-for-the-mean} is called \emph{Pettis integral}; moreover, in general, there exists a 
continuous linear operator $\Sigma_X:B^*\to B^{**}$, where $B^{**}$ is the dual space of $B^*$ (however, in our setting, $\Sigma_X$ take 
values in $B$).

In this paper we consider the following Condition \ref{cond:LD-*} which plays a crucial role in view of the proof of the 
reference LDP (Proposition \ref{prop:LD}). Such a condition allows to consider a straightforward application of the 
G\"artner--Ellis Theorem (see, e.g., Theorem 2.3.6 in \cite{DZ:97}) for $\{N_n/n:n\geq 1\}$; so, in some sense, Condition \ref{cond:LD-*}
can be seen as the analogue of Assumption 4.1 in \cite{EicheLo}.

\begin{condition}\label{cond:LD-*}
Let $\{N_n:n\geq 1\}$ be a family of integer-valued nonnegative random variables such that the function
$$\Lambda_N(\eta):=\lim_{n\to\infty}\frac{1}{n}\log\mathbb{E}[e^{\eta N_n}]\quad\quad\mbox{for all}\ \eta\in\mathbb{R}$$
is well-defined, finite valued and differentiable; note that, since $\Lambda_N$ is a non-decreasing function,
it is well-defined $\Lambda_N(-\infty):=\lim_{\eta\to-\infty}\Lambda_N(\eta)$. Moreover, let $\{X_n:n\geq 1\}$ be a 
sequence of i.i.d. $B$-valued random variables as in Condition \ref{cond:superexponential-tails-summands}, and 
independent of $\{N_n:n\geq 1\}$.
\end{condition}

Moreover we consider the assumptions collected in the next Condition \ref{cond:MD-*}. These assumptions play a crucial role
in the proofs of the moderate deviation results (Propositions \ref{prop:MD-summands-are-centered} and 
\ref{prop:MD-centered-compound-sum}); in particular the second one with eq. \eqref{eq:MD-limit-for-marginal-N}
(which, in some sense, reminds a condition in Laplace's method) allows to have a moderate deviation result for the random 
variables $\{N_n:n\geq 1\}$.

\begin{condition}\label{cond:MD-*}
Assume that Condition \ref{cond:LD-*} holds, and there exists $\Lambda_N^{\prime\prime}(0)$. Moreover assume that
\begin{itemize}
\item $\frac{S_n-n\mu_X}{\sqrt{n}}$ converges weakly to a centered normal distribution (as $n\to\infty$);
\item given a positive sequence $\{a_n:n\geq 1\}$ such that eq. \eqref{eq:MD-conditions} holds, we have
\begin{equation}\label{eq:MD-limit-for-marginal-N}
\lim_{n\to\infty}a_n\log\mathbb{E}\left[e^{\eta\frac{N_n-\mathbb{E}[N_n]}{\sqrt{na_n}}}\right]
=\frac{\Lambda_N^{\prime\prime}(0)}{2}\eta^2\quad\quad\mbox{for all}\ \eta\in\mathbb{R}.
\end{equation}
\end{itemize}
\end{condition}
We remark that, as a consequence of Condition \ref{cond:MD-*}, we have exponential tightness results for the
moderate deviations of $\{S_n:n\geq 1\}$ (see \cite{deAcosta3}) and of $\{N_n:n\geq 1\}$; those results will be used to prove 
the exponential tightness for the moderate deviation results in this paper. We also have an analogue remark for Condition 
\ref{cond:LD-*} and the exponential tightness for the reference LDP in this paper (in this case the bibliographic reference 
for $\{S_n:n\geq 1\}$ is \cite{deAcosta1}).

We can also say that $\Lambda_N^{\prime\prime}(0)\geq 0$. However, in order to avoid 
the degenerating cases in which the random sequence $\{N_n:n\geq 1\}$ behaves as a deterministic sequence, we require
that $\Lambda_N^{\prime\prime}(0)>0$; see, e.g., Proposition \ref{prop:MD-summands-are-centered}, in which the rate 
function expression has a denominator equal to $\Lambda_N^{\prime\prime}(0)$ (obviously this also happens for the
rate function expression in Proposition \ref{prop:MD-centered-compound-sum}).

For completeness we remark that the class of all positive scalings $\{a_n:n\geq 1\}$ (such that eq. 
\eqref{eq:MD-conditions} holds) allows to fill the gap between two asymptotic regimes as $n\to\infty$: 
\begin{enumerate}
	\item the convergence of $\left(\frac{S_{N_n}}{n},\frac{N_n}{n}\right)$ to $(\Lambda_N^\prime(0)\mu_X,\Lambda_N^\prime(0))$ (at 
	least in probability) governed by a LDP (i.e. the \emph{reference LDP} cited above);
	\item the weak convergence of $\left(\frac{S_{N_n}-N_n\mu_X}{\sqrt{n}},\frac{N_n-\mathbb{E}[N_n]}{\sqrt{n}}\right)$ and $\left(\frac{S_{N_n}-\mathbb{E}[N_n]\mu_X}{\sqrt{n}},\frac{N_n-\mathbb{E}[N_n]}{\sqrt{n}}\right)$ to two centered normal distributions
	on $B\times\mathbb{R}$ (see also Remark \ref{rem:on-the-weak-convergences} for more details).
\end{enumerate}
Moreover, we recover the first and the second regime by taking $a_n=\frac{1}{v_n}=\frac{1}{n}$ and $a_n=1$, respectively; note that,
in these two cases, one condition in eq. \eqref{eq:MD-conditions} holds, and the other one fails.

\section{Large deviations}\label{sec:LD}
In this section we prove Proposition \ref{prop:LD}. In particular we show that we have an exponentially tight sequence, 
and this will be done separately in Section \ref{sub:LD-ET}.

\subsection{The result and a part of the proof}
We have the following result.

\begin{proposition}\label{prop:LD}
Assume that Condition \ref{cond:LD-*} holds. Then the sequence
$\left\{\left(\frac{S_{N_n}}{n},\frac{N_n}{n}\right):n\geq 1\right\}$ satisfies the LDP 
with speed $n$, and good rate function $I_{S,N}$ defined by
\begin{equation}\label{eq:LD-rf-LT}
I_{S,N}(x,y):=\sup_{\theta\in B^*,\eta\in\mathbb{R}}\{\langle\theta,x\rangle+\eta y-\Lambda_N(\eta+\Lambda_X(\theta))\}
\quad\quad\mbox{for all}\ (x,y)\in B\times\mathbb{R}.
\end{equation}
Moreover, we have
\begin{equation}\label{eq:LD-rf-explicit}
I_{S,N}(x,y):=\left\{\begin{array}{ll}
y\Lambda_X^*\left(\frac{x}{y}\right)+\Lambda_N^*(y)&\ \mbox{if}\ y>0\\
-\Lambda_N(-\infty)&\ \mbox{if}\ (x,y)=(\widehat{0},0)\\
\infty&\ \mbox{otherwise}.
\end{array}\right.
\end{equation}
\end{proposition}
\begin{proof}
We prove the LDP as an application of Corollary 4.6.14 in \cite{DZ:97} (it is a stronger version of Corollary 4.5.27 in 
\cite{DZ:97}, which is a consequence of Theorem 4.5.20 in the same reference). Firstly, by the G\"artner--Ellis Theorem, 
$\left\{\frac{N_n}{n}:n\geq 1\right\}$ satisfies the LDP with good rate function $\Lambda_N^*$ defined by
$$\Lambda_N^*(y):=\sup_{\eta\in\mathbb{R}}\{\eta y-\Lambda_N(\eta)\}.$$
Then $\{\frac{N_n}{n}:n\geq 1\}$ is exponentially tight at speed $n$ (see, for example, Lemma 2.6 in \cite{Lyn-Set}); therefore,
for all $R>0$, there exists $L_R>0$ such that, for some $n_0\geq 1$,
\begin{equation}\label{eq:LD-ET-counting-processes}
P\left(\frac{N_n}{n}>L_R\right)\leq e^{-nR}\quad\quad\mbox{for all}\ n\geq n_0.
\end{equation}
Moreover, as we shall explain in Subsection \ref{sub:LD-ET}, the sequence
$\left\{\left(\frac{S_{N_n}}{n},\frac{N_n}{n}\right):n\geq 1\right\}$ is exponentially tight at speed $n$
(this is stated in Proposition \ref{prop:LD-ET}).

Now we remark that
$$\mathbb{E}[e^{\langle\theta,S_{N_n}\rangle+\eta N_n}]
=\mathbb{E}[\mathbb{E}[e^{\langle\theta,X_1\rangle}]^{N_n}e^{\eta N_n}]
=\mathbb{E}[e^{(\eta+\Lambda_X(\theta))N_n}]\quad\quad\mbox{for all}\ (\theta,\eta)\in B^*\times\mathbb{R};$$
thus we have
\begin{equation}\label{eq:def-LambdaSN}
\Lambda_{S,N}(\theta,\eta):=\lim_{n\to\infty}\frac{1}{n}\log\mathbb{E}[e^{\langle\theta,S_{N_n}\rangle+\eta N_n}]
=\Lambda_N(\eta+\Lambda_X(\theta))\quad\quad\mbox{for all}\ (\theta,\eta)\in B^*\times\mathbb{R}.
\end{equation}
Then, by taking into account our hypotheses, the function $\Lambda_{S,N}$ is obviously finite valued and Gateaux 
differentiable; thus, by Corollary 4.6.14 in \cite{DZ:97}, we get the desired LDP with good rate function 
$I_{S,N}$ defined by eq. \eqref{eq:LD-rf-LT}.

In the final part of the proof we show that eq. \eqref{eq:LD-rf-explicit} holds. We have the following cases.
\begin{itemize}
\item If $y<0$, then
$$I_{S,N}(x,y)\geq\langle\widehat{0},x\rangle+\eta y-\Lambda_N(\eta+\Lambda_X(\widehat{0}))
=\eta y-\Lambda_N(\eta);$$
moreover, since $\Lambda_N(\eta)\leq 0$ for $\eta\leq 0$ (because $\Lambda_N$ is a non-decreasing function),
for $\eta\leq 0$ we have $I_{S,N}(x,y)\geq\eta y$ and we conclude by letting $\eta$ go to $-\infty$.
\item If $y>0$, then
\begin{multline*}
I_{S,N}(x,y)\leq\sup_{\theta\in B^*}\{\langle\theta,x\rangle-y\Lambda_X(\theta)\}
+\sup_{\theta\in B^*,\eta\in\mathbb{R}}\{(\eta+\Lambda_X(\theta))y-\Lambda_N(\eta+\Lambda_X(\theta))\}\\
=y\sup_{\theta\in B^*}\{\langle\theta,x/y\rangle-\Lambda_X(\theta)\}+\Lambda_N^*(y)
=y\Lambda_X^*\left(\frac{x}{y}\right)+\Lambda_N^*(y);
\end{multline*}
moreover, if we take two maximizing sequences $\{\theta_n:n\geq 1\}$ and $\{\alpha_n:n\geq 1\}$ such that
$$\Lambda_X^*\left(\frac{x}{y}\right)=\lim_{n\to\infty}\langle\theta_n,x/y\rangle-\Lambda_X(\theta_n)\quad
\mbox{and}\quad \Lambda_N^*(y)=\lim_{n\to\infty}\alpha_ny-\Lambda_N(\alpha_n),$$
for every $n\geq 1$ we get
\begin{multline*}
I_{S,N}(x,y)\geq\left.\langle\theta_n,x\rangle+\eta y-\Lambda_N(\eta+\Lambda_X(\theta_n))\right|_{\eta=-\Lambda_X(\theta_n)+\alpha_n}\\
=\langle\theta_n,x\rangle-y\Lambda_X(\theta_n)+\alpha_ny-\Lambda_N(\alpha_n)
=y(\langle\theta_n,x/y\rangle-\Lambda_X(\theta_n))+\alpha_ny-\Lambda_N(\alpha_n),
\end{multline*}
and we conclude by letting $n$ go to infinity (because we have shown that the inverse ineqality holds).
\item If $y=0$ and $x=\widehat{0}$, then
$$I_{S,N}(\widehat{0},0)=\sup_{\theta\in B^*,\eta\in\mathbb{R}}\{-\Lambda_N(\eta+\Lambda_X(\theta))\}
=-\lim_{\eta\to-\infty}\Lambda_N(\eta)=-\Lambda_N(-\infty)$$
because $\Lambda_N$ is a non-decreasing function.
\item If $y=0$ and $x\neq\widehat{0}$, then we take $\theta_x\in B^*$ such that $\langle\theta_x,x\rangle>0$;
thus for every $n\geq 1$ we get
$$I_{S,N}(\widehat{0},0)\geq
\left.\langle n\theta_x,x\rangle+\eta\cdot 0-\Lambda_N(\eta+\Lambda_X(n\theta_x))\right|_{\eta=-\Lambda_X(n\theta_x)}
=n\langle\theta_x,x\rangle-\Lambda_N(0)=n\langle\theta_x,x\rangle,$$
and we conclude by letting $n$ go to infinity.
\end{itemize}
\end{proof}

\begin{remark}\label{rem:typical-features-for-N-only}
By the hypotheses and Theorem 25.7 in \cite{Rockafeller} we have
\begin{equation}\label{eq:limit-mean-tf-N}
\lim_{n\to\infty}\mathbb{E}[N_n/n]=\Lambda_N^\prime(0).
\end{equation}
Moreover, one typically has
\begin{equation}\label{eq:limit-variance-tf-N}
\lim_{n\to\infty}n\mathrm{Var}[N_n/n]=\Lambda_N^{\prime\prime}(0).
\end{equation}
In Section \ref{sub:some-typical-features} we present suitable versions of the limits in eqs. \eqref{eq:limit-mean-tf-N} and 
\eqref{eq:limit-variance-tf-N} for the bivariate sequence in Proposition \ref{prop:LD}.
\end{remark}

\subsection{The exponential tightness of $\left\{\left(\frac{S_{N_n}}{n},\frac{N_n}{n}\right):n\geq 1\right\}$ at speed $n$}\label{sub:LD-ET}
We start with some preliminaries. For some $K\subset B$, we introduce the Minkowski functional
$$q_K(x):=\inf\{a>0:x\in aK\},$$
where, as usual, $\inf\emptyset=\infty$. We say that a set $K$ is well-balanced if $aK\subset K$ for $a\in[0,1]$. 
If a set is well-balanced then $a_1K\subset a_2K$ for $0<a_1\leq a_2$; moreover $x\notin aK$ yields $q_K(x)>a$. 
Furthermore recall that for a convex well-balanced set, the Minkowski functional is a semi-norm. See, e.g., 
\cite{Schechter} or \cite{Nar-Beck}.

Then we have the following lemmas.

\begin{lemma}\label{lem:LD-*}
Assume that Condition \ref{cond:superexponential-tails-summands} holds. Then there exists a well-balanced 
compact convex set $K\subset B$ such that $\mathbb{E}[e^{q_K(X_1)}]<\infty$.
\end{lemma}
\begin{proof}
See, e.g., Theorem 3.1 in \cite{deAcosta1}.
\end{proof}

\begin{lemma}\label{lem:LD-**}
Assume that Condition \ref{cond:superexponential-tails-summands} holds, and let $K$ be the set in Lemma \ref{lem:LD-*}.
Then, for $R,L>0$, there exists $\delta=\delta(R,L,K)>0$ such that
$$\sup_{1\leq h\leq Ln}P\left(\frac{S_h}{n}\notin\delta K\right)\leq e^{-nR}\quad\quad\mbox{for all integers}\ n\geq 1.$$
\end{lemma}
\begin{proof}
Let $\delta>0$ be arbitrarily fixed. Then, for $n\geq 1$ and for all integers $h\in[1,Ln]$, we have
\begin{multline*}
P\left(\frac{S_h}{n}\notin\delta K\right)=P(S_h\notin n\delta K)=P(q_K(S_h)>\delta n)\\
\leq\mathbb{E}[e^{q_K(S_h)}]e^{-n\delta}\leq\mathbb{E}^h[e^{q_K(X_1)}]e^{-n\delta}
=e^{-n(\delta-\frac{h}{n}\log\mathbb{E}[e^{q_K(X_1)}])};
\end{multline*}
here we take into account that $K$ is a well-balanced set and the sublinearity of $q_K(x)$ 
because $q_K$ is a semi-norm. So we conclude by taking $\delta=\delta(R,L,K)>0$ such that
$$\delta-\frac{h}{n}\log\mathbb{E}[e^{q_K(X_1)}]>R\quad\quad\mbox{for all integers}\ h\in[1,Ln];$$
indeed this is possible by taking $\delta>R+L\max\{\log\mathbb{E}[e^{q_K(X_1)}],0\}$.
\end{proof}

Now we are able to prove the desired exponential tightness condition.

\begin{proposition}\label{prop:LD-ET}
Assume that Condition \ref{cond:LD-*} holds. Then the sequence 
$\left\{\left(\frac{S_{N_n}}{n},\frac{N_n}{n}\right):n\geq 1\right\}$ 
is exponentially tight at speed $n$.
\end{proposition}
\begin{proof}
In this proof we show that, for all $R>0$, there exists a compact set $\widehat{K}_R\subset B\times\mathbb{R}$
such that, for some $n_2\geq 1$,
$$P\left(\left(\frac{S_{N_n}}{n},\frac{N_n}{n}\right)\notin\widehat{K}_R\right)\leq 2e^{-nR}
\quad\quad\mbox{for all}\ n\geq n_2.$$	
More precisely we consider the set
$$\widehat{K}_R:=(\delta_R K)\times [0,L_R],$$
where $\delta_R:=\delta(R,L_R,K)>0$ according to the notation in Lemma \ref{lem:LD-**} (with the compact $K$ in that lemma),
and $L_R>0$ as in eq. \eqref{eq:LD-ET-counting-processes}. Then we have
\begin{multline*}
P\left(\left(\frac{S_{N_n}}{n},\frac{N_n}{n}\right)\notin\widehat{K}_R\right)
=P\left(\left(\left\{\frac{S_{N_n}}{n}\in\delta_RK\right\}\cap
\left\{\frac{N_n}{n}\in[0,L_R]\right\}\right)^c\right)\\
=P\left(\left\{\frac{S_{N_n}}{n}\notin\delta_RK\right\}\cup
\left\{\frac{N_n}{n}>L_R\right\}\right)\\
=P\left(\frac{N_n}{n}>L_R\right)+P\left(\left\{\frac{S_{N_n}}{n}\notin\delta_RK\right\}\cap
\left\{\frac{N_n}{n}\leq L_R\right\}\right).
\end{multline*}
Now we estimate the summands in the final expression. We refer to eq. \eqref{eq:LD-ET-counting-processes}
to estimate the first summand. For the second one we take into account Lemma \ref{lem:LD-**} with $L=L_R$
(and the definition of $\delta_R$ in this proof); then, for some $n_1\geq 1$, we have (here, without loss of generality, $L_R$
is integer)
\begin{multline*}
P\left(\left\{\frac{S_{N_n}}{n}\notin\delta_RK\right\}\cap
\left\{\frac{N_n}{n}\leq L_R\right\}\right)
=\sum_{j=0}^{nL_R}P\left(\left\{\frac{S_{N_n}}{n}\notin\delta_RK\right\}\cap\{N_n=j\}\right)\\
=\sum_{j=0}^{nL_R}P\left(\frac{S_j}{n}\notin\delta_RK\right)P(N_n=j)
\leq e^{-nR}\sum_{j=0}^{nL_R}P(N_n=j)\leq e^{-nR}.
\end{multline*}
Thus, if we take $n_2:=\max\{n_0,n_1\}$, the statement at the beginning of the proof is proved.
This completes the proof. 
\end{proof}

\section{Moderate deviations (MD)}\label{sec:MD}
In this Section we prove two results: Proposition \ref{prop:MD-summands-are-centered} (Sections 
\ref{sub:MD-summands-are-centered} and \ref{sub:MD-summands-are-centered-ET}) for the case in which ``the summands 
are centered'' and Proposition \ref{prop:MD-centered-compound-sum} (Section \ref{sub:MD-centered-compound-sum})
for the case in which ``the compound sum is centered''.
In the first case, for $\mu_X\in B$ as in eq. \eqref{eq:condition-for-the-mean}, we mean
$$\sum_{i=1}^{N_n}(X_i-\mu_X)=S_{N_n}-N_n\mu_X;$$
in the second case we mean
$$S_{N_n}-\mathbb{E}[S_{N_n}]=S_{N_n}-\mathbb{E}[N_n]\mu_X;$$
obviously these two cases coincide if and only if $\mu_X=\widehat{0}$. In particular, for the proof of Proposition 
\ref{prop:MD-summands-are-centered}, in Section \ref{sub:MD-summands-are-centered-ET} we prove an exponential tightness condition.
Proposition \ref{prop:MD-centered-compound-sum} will be a consequence of Proposition \ref{prop:MD-summands-are-centered} and an 
application of the contraction principle. A final Section \ref{sub:some-typical-features} is devoted to the discussion of some 
typical features.

Finally we recall that, if $\langle\theta,\Sigma_X\theta\rangle>0$ for all $\theta\neq\widehat{0}^*$, then
we can define the inverse of the function $\Sigma_X$ in eq. \eqref{eq:condition-for-the-covariance}; more precisely we mean 
$\Sigma_X^{-1}:\mathrm{Im}(\Sigma_X)\to B^*$, where
$$\mathrm{Im}(\Sigma_X):=\{x\in B:\Sigma_Xu=x\ \mbox{for some}\ u\in B^*\},$$
such that $\Sigma_X^{-1}\Sigma_Xu=u$ for all $u\in B^*$.

\subsection{MD when ``the summands are centered'' and a part of the proof}\label{sub:MD-summands-are-centered}
We have the following result.

\begin{proposition}\label{prop:MD-summands-are-centered}
Let $\mu_X\in B$ and the linear operator $\Sigma_X$ be defined in eqs. \eqref{eq:condition-for-the-mean} and 
\eqref{eq:condition-for-the-covariance}, respectively. Assume that Condition \ref{cond:MD-*} and the limit in eq. 
\eqref{eq:limit-mean-tf-N} hold.
Then, for every positive sequence $\{a_n:n\geq 1\}$ such that eq. \eqref{eq:MD-conditions} holds, the sequence 
$\left\{\left(\frac{S_{N_n}-N_n\mu_X}{\sqrt{n/a_n}},\frac{N_n-\mathbb{E}[N_n]}{\sqrt{n/a_n}}\right):n\geq 1\right\}$ satisfies
the LDP with speed $1/a_n$, and good rate function $J_{S,N}^{(1)}$ defined by
$$J_{S,N}^{(1)}(x,y):=\sup_{\theta\in B^*}\left\{\langle\theta,x\rangle-\frac{\Lambda_N^\prime(0)\langle\theta,\Sigma_X\theta\rangle}{2}\right\}
+\frac{y^2}{2\Lambda_N^{\prime\prime}(0)}.$$
Moreover we have the following particular cases:
$$J_{S,N}^{(1)}(x,y):=\left\{\begin{array}{ll}
\frac{y^2}{2\Lambda_N^{\prime\prime}(0)}&\ \mbox{if}\ x=\widehat{0}\\
\infty&\ \mbox{otherwise}
\end{array}\right.$$
if $\Lambda_N^\prime(0)=0$;
$$J_{S,N}^{(1)}(x,y):=\left\{\begin{array}{ll}
	\frac{\langle x,\Sigma_X^{-1}x\rangle}{2\Lambda_N^\prime(0)}
	+\frac{y^2}{2\Lambda_N^{\prime\prime}(0)}&\ \mbox{if}\ x\in\mathrm{Im}(\Sigma_X)\\
	\infty&\ \mbox{otherwise}
\end{array}\right.$$
if $\Lambda_N^\prime(0)>0$ and $\langle\theta,\Sigma_X\theta\rangle>0$ for all $\theta\neq\widehat{0}^*$.
\end{proposition}
\begin{proof}
We start with the particular cases in the final part of the statement. The first one is immediate, while the second one holds noting 
that $\sup_{\theta\in B^*}\left\{\langle\theta,x\rangle-\frac{\Lambda_N^\prime(0)\langle\theta,\Sigma_X\theta\rangle}{2}\right\}$
is attained at $\theta=\frac{1}{2\Lambda_N^\prime(0)}\langle\Sigma_X^{-1}x,x\rangle$ (and some easy computations).

We consider an arbitrarily fixed positive sequence $\{a_n:n\geq 1\}$ such that eq. \eqref{eq:MD-conditions} holds and, as we did 
for Proposition \ref{prop:LD}, we prove the LDP as an application of Corollary 4.6.14 in \cite{DZ:97}. Firstly, 
by the G\"artner--Ellis Theorem and the limit in eq. \eqref{eq:MD-limit-for-marginal-N}, 
$\left\{\frac{N_n-\mathbb{E}[N_n]}{\sqrt{n/a_n}}:n\geq 1\right\}$ satisfies the LDP with good rate function $J_N$ defined by
$$J_N^*(y):=\sup_{\eta\in\mathbb{R}}\left\{\eta y-\frac{\Lambda_N^{\prime\prime}(0)}{2}\eta^2\right\}
=\frac{y^2}{2\Lambda_N^{\prime\prime}(0)}.$$
Then, again by Lemma 2.6 in \cite{Lyn-Set}, $\{\frac{N_n-\mathbb{E}[N_n]}{\sqrt{n/a_n}}:n\geq 1\}$ is exponentially 
tight at speed $1/a_n$; therefore, for all $R>0$, there exists $L_R>0$ such that, for some $n_0\geq 1$,
\begin{equation}\label{eq:MD-summands-are-centered-ET-counting-processes}
P\left(\frac{|N_n-\mathbb{E}[N_n]|}{\sqrt{n/a_n}}>L_R\right)\leq e^{-R/a_n}\quad\quad\mbox{for all}\ n\geq n_0.
\end{equation}
Moreover, as we shall explain in Subsection \ref{sub:MD-summands-are-centered-ET}, the sequence
$\left\{\left(\frac{S_{N_n}-N_n\mu_X}{\sqrt{n/a_n}},\frac{N_n-\mathbb{E}[N_n]}{\sqrt{n/a_n}}\right):n\geq 1\right\}$ is 
exponentially tight at speed $1/a_n$ (this is stated in Proposition \ref{prop:MD-summands-are-centered-ET}).

Now we remark
\begin{multline*}
\mathbb{E}[e^{\langle\theta,S_{N_n}-N_n\mu_X\rangle+\eta (N_n-\mathbb{E}[N_n])}]\\
=\mathbb{E}[\mathbb{E}[e^{\langle\theta,X_1-\mu_X\rangle}]^{N_n}e^{\eta N_n}]e^{-\eta\mathbb{E}[N_n]}
=\mathbb{E}[e^{(\eta+\Lambda_X(\theta)-\langle\theta,\mu_X\rangle)N_n}]e^{-\eta\mathbb{E}[N_n]}\\
=\mathbb{E}[e^{(\eta+\Lambda_X(\theta)-\langle\theta,\mu_X\rangle)(N_n-\mathbb{E}[N_n])}]
e^{(\Lambda_X(\theta)-\langle\theta,\mu_X\rangle)\mathbb{E}[N_n]}\quad\quad\mbox{for all}\ (\theta,\eta)\in B^*\times\mathbb{R};
\end{multline*}
thus we have
\begin{multline*}
a_n\log\mathbb{E}\left[e^{\langle\theta,\frac{S_{N_n}-N_n\mu_X}{\sqrt{na_n}}\rangle+\eta\frac{N_n-\mathbb{E}[N_n]}{\sqrt{na_n}}}\right]\\
=a_n\log\mathbb{E}[e^{(\eta/\sqrt{na_n}+\Lambda_X(\theta/\sqrt{na_n})-\langle\theta/\sqrt{na_n},\mu_X\rangle)(N_n-\mathbb{E}[N_n])}]
+a_n(\Lambda_X(\theta/\sqrt{na_n})-\langle\theta/\sqrt{na_n},\mu_X\rangle)\mathbb{E}[N_n].
\end{multline*}

Then, in what follows, we estimate the summands in the final expression. Actually it is easy to see 
(and we omit the details) that for the second term we have
$$\lim_{n\to\infty}a_n(\Lambda_X(\theta/\sqrt{na_n})-\langle\theta/\sqrt{na_n},\mu_X\rangle)\mathbb{E}[N_n]
=\frac{\Lambda_N^\prime(0)}{2}\langle\theta,\Sigma_X\theta\rangle\quad\quad\mbox{for all}\ \theta\in B^*,$$
where we take into account the limit in eq. \eqref{eq:limit-mean-tf-N}. So we complete the proof showing that
\begin{equation}\label{eq:MD-summands-are-centered-limit-to-prove}
\lim_{n\to\infty}a_n\log\mathbb{E}[e^{(\eta/\sqrt{na_n}+\Lambda_X(\theta/\sqrt{na_n})-\langle\theta/\sqrt{na_n},\mu_X\rangle)(N_n-\mathbb{E}[N_n])}]
=\frac{\Lambda_N^{\prime\prime}(0)}{2}\eta^2\quad\quad\mbox{for all}\ (\theta,\eta)\in B^*\times\mathbb{R};
\end{equation}
indeed, if we consider the function $\Psi_{S,N}$ defined by
\begin{equation}\label{eq:def-PsiSN}
\Psi_{S,N}(\theta,\eta):=\frac{\Lambda_N^\prime(0)\langle\theta,\Sigma_X\theta\rangle}{2}+\frac{\Lambda_N^{\prime\prime}(0)}{2}\eta^2
\end{equation}
(which is obviously finite valued and Gateaux differentiable), the desired LDP holds by Corollary 4.6.14 in \cite{DZ:97} with 
good rate function $J_{S,N}^{(1)}$ defined by
\begin{equation}\label{eq:MD-rf-LT-summands-are-centered}
J_{S,N}^{(1)}(x,y):=\sup_{\theta\in B^*,\eta\in\mathbb{R}}\{\langle\theta,x\rangle+\eta y-\Psi_{S,N}(\theta,\eta)\}
\quad\quad\mbox{for all}\ (x,y)\in B\times\mathbb{R},
\end{equation}
which coincides with the rate function $J_{S,N}^{(1)}$ because
$$\sup_{\theta\in B^*,\eta\in\mathbb{R}}\{\langle\theta,x\rangle+\eta y-\Psi_{S,N}(\theta,\eta)\}
=\sup_{\theta\in B^*}\left\{\langle\theta,x\rangle-\frac{\Lambda_N^\prime(0)\langle\theta,\Sigma_X\theta\rangle}{2}\right\}
+\sup_{\eta\in\mathbb{R}}\left\{\eta y-\frac{\Lambda_N^{\prime\prime}(0)}{2}\eta^2\right\}.$$
Thus we complete the proof showing that \eqref{eq:MD-summands-are-centered-limit-to-prove} holds. The case
$\theta=\widehat{0}^*$ is trivial because it is an immediate consequence of the limit in eq. \eqref{eq:MD-limit-for-marginal-N}.
So, from now on, we can restrict the attention on the case $\theta\neq\widehat{0}^*$. Moreover we note that, if we set
$$\eta_n=\eta+\frac{1}{2\sqrt{na_n}}\langle\theta,\Sigma_X\theta\rangle+\sqrt{na_n}o\left(\frac{1}{na_n}\right),$$
we have $\lim_{n\to\infty}\eta_n=\eta$ and
$$\mathbb{E}[e^{(\eta/\sqrt{na_n}+\Lambda_X(\theta/\sqrt{na_n})-\langle\theta/\sqrt{na_n},\mu_X\rangle)(N_n-\mathbb{E}[N_n])}]
=\mathbb{E}[e^{\eta_n/\sqrt{na_n}(N_n-\mathbb{E}[N_n])}];$$
so we get eq. \eqref{eq:MD-summands-are-centered-limit-to-prove} (and we complete the proof) if we show that
\begin{equation}\label{eq:MD-summands-are-centered-limit-to-prove-2}
\lim_{n\to\infty}a_n\log\mathbb{E}[e^{\eta_n/\sqrt{na_n}(N_n-\mathbb{E}[N_n])}]
=\frac{\Lambda_N^{\prime\prime}(0)}{2}\eta^2\quad\quad\mbox{for all}\ (\theta,\eta)\in (B^*\setminus\{\widehat{0}^*\})\times\mathbb{R}.
\end{equation}
This desired limit in eq. \eqref{eq:MD-summands-are-centered-limit-to-prove-2} can be proved noting that 
the convergence in eq. \eqref{eq:MD-limit-for-marginal-N} holds uniformly on compact sets (by Theorem 10.8 in \cite{Rockafeller})
and some standard computations.
\end{proof}

\subsection{The exponential tightness of 
$\left\{\left(\frac{S_{N_n}-N_n\mu_X}{\sqrt{n/a_n}},\frac{N_n-\mathbb{E}[N_n]}{\sqrt{n/a_n}}\right):n\geq 1\right\}$ 
at speed $1/a_n$}\label{sub:MD-summands-are-centered-ET}
Here we refer again to the Minkowski functional $q_K$, for some $K\subset B$, and some of its properties (see the 
beginning of Section \ref{sub:LD-ET}). We have the following lemmas and, as usual, we refer to $\mu_X\in B$ defined
in eq. \eqref{eq:condition-for-the-mean}.

\begin{lemma}\label{lem:MD-summands-are-centered-*}
Assume that Condition \ref{cond:MD-*} holds. Let $\{a_n:n\geq 1\}$ be a positive sequence such that eq. 
\eqref{eq:MD-conditions} holds. Then there exists a well-balanced compact convex set $K\subset B$ such that, for every $L>0$, 
there exists a constant $c=c(L,K)>0$ such that
$$\sup_{1\leq h\leq Ln}\mathbb{E}[e^{q_K(S_h-h\mu_X)/\sqrt{na_n}}]\leq c^{1/a_n}.$$
\end{lemma}
\begin{proof}
We start noting that, by taking into account the first part	of the proof of Lemma 2.1 in \cite{deAcosta3} 
(we recall that $B$ is separable, and therefore the law of $X_1$ is tight as explained in Condition 
\ref{cond:superexponential-tails-summands}), we can say what follows:
$$\beta:=\sup\left\{\mathbb{E}\left[\exp\left(\frac{q_K(S_n-n\mu_X)}{\sqrt{n}}\right)\right]:n\geq 1\right\}<\infty;$$
for every $\lambda\in(0,1)$ we have
$$\mathbb{E}\left[\exp\left(\lambda(q_K(S_n-n\mu_X)-\mathbb{E}[q_K(S_n-n\mu_X)])\right)\right]\leq\exp(\lambda^2n\alpha)
\quad\quad\mbox{for all integers}\ n\geq 1,$$
where
$$\alpha:=\mathbb{E}[\exp(2q_K(X_1-\mu_X)-\mathbb{E}[q_K(X_1-n\mu_X)])]<\infty$$
(indeed $\alpha$ is finite because $\mathbb{E}[\exp(2q_K(X_1-\mu_X))]<\infty$, as shown in the proof of Lemma 2.1 in 
\cite{deAcosta3}). Then, for $n$ large enough, we take $\lambda=\frac{1}{\sqrt{na_n}}$ and, for all integers $h\in[1,Ln]$,
we have
$$\mathbb{E}\left[\exp\left(\frac{1}{\sqrt{na_n}}(q_K(S_h-h\mu_X)-\mathbb{E}[q_K(S_h-h\mu_X)])\right)\right]
\leq\exp(h\alpha/(na_n))\leq e^{L\alpha/a_n};$$
therefore (here we take into account the sublinearity of $q_K(x)$ because $q_K$ is a semi-norm)
$$\mathbb{E}[\exp(q_K(S_h-h\mu_X)/\sqrt{na_n})]\leq\mathbb{E}[\exp(\mathbb{E}[q_K((S_h-h\mu_X)/\sqrt{n})]/\sqrt{a_n})]e^{L\alpha/a_n}.$$
Moreover, again for all integers $h\in[1,Ln]$, for  $\beta$ defined above we get
$$\mathbb{E}[q_K((S_h-h\mu_X)/\sqrt{n})]=\mathbb{E}[\sqrt{h/n}q_K((S_h-h\mu_X)/\sqrt{h})]
\leq\sqrt{L}\mathbb{E}[q_K((S_h-h\mu_X)/\sqrt{h})]\leq\sqrt{L}\beta;$$
therefore
$$\mathbb{E}[e^{q_K(S_h-h\mu_X)/\sqrt{na_n}}]\leq e^{\sqrt{L}\beta/\sqrt{a_n}+L\alpha/a_n}\leq e^{(\sqrt{L}\beta+L\alpha)/a_n},$$
and the proof is complete by setting $c=e^{\sqrt{L}\beta+L\alpha}$.
\end{proof}

\begin{lemma}\label{lem:MD-summands-are-centered-**}
Assume that Condition \ref{cond:MD-*} holds, and let $K$ be the set in Lemma \ref{lem:MD-summands-are-centered-*}.
Let $\{a_n:n\geq 1\}$ be a positive sequence such that eq. \eqref{eq:MD-conditions} holds.
Then, for $R,L>0$, there exists $\delta=\delta(R,L,K)>0$ such that
$$\sup_{1\leq h\leq Ln}P\left(\frac{S_h-h\mu_X}{\sqrt{n/a_n}}\notin\delta K\right)\leq e^{-R/a_n}\quad\quad\mbox{for all integers}\ n\geq 1.$$
\end{lemma}
\begin{proof}
Let $c$ be the quantity in Lemma \ref{lem:MD-summands-are-centered-*} (which depends on $L$) and let $\delta>0$ be arbitrarily 
fixed. Then, for $n\geq 1$ and for all integers $h\in[1,Ln]$, we have
\begin{multline*}
P\left(\frac{S_h-h\mu_X}{\sqrt{n/a_n}}\notin\delta K\right)=P\left(\frac{S_h-h\mu_X}{\sqrt{na_n}}\notin\frac{\delta}{a_n}K\right)
=P\left(\frac{q_K(S_h-h\mu_X)}{\sqrt{na_n}}>\frac{\delta}{a_n}\right)\\
\leq\mathbb{E}[e^{q_K(S_h-h\mu_X)/\sqrt{na_n}}]e^{-\delta/a_n}\leq c^{1/a_n}e^{-\delta/a_n}
=e^{-(\delta-\log c)/a_n};
\end{multline*}
(in the second equality we take into account that $K$ is a well-balanced set and the sublinearity of $q_K(x)$ 
because $q_K$ is a semi-norm; then we consider the Markov inequality and Lemma \ref{lem:MD-summands-are-centered-*}). 
So we conclude by taking $\delta=\delta(R,L,K)>0$ such that $\delta-\log c>R$.
\end{proof}

Now we are able to prove the desired exponential tightness condition.

\begin{proposition}\label{prop:MD-summands-are-centered-ET}
Assume that Condition \ref{cond:MD-*} and the limit in eq. \eqref{eq:limit-mean-tf-N} hold.
Let $\{a_n:n\geq 1\}$ be a positive sequence such that eq. \eqref{eq:MD-conditions} holds. Then 
the sequence $\left\{\left(\frac{S_{N_n}}{n},\frac{N_n}{n}\right):n\geq 1\right\}$ is exponentially tight at speed $1/a_n$.
\end{proposition}
\begin{proof}
In this proof we show that, for all $R>0$, there exists a compact set $\widehat{K}_R\subset B\times\mathbb{R}$
such that, for some $n_2\geq 1$,
$$P\left(\left(\frac{S_{N_n}-N_n\mu_X}{\sqrt{n/a_n}},\frac{N_n-\mathbb{E}[N_n]}{\sqrt{n/a_n}}\right)\notin\widehat{K}_R\right)\leq 2e^{-R/a_n}
\quad\quad\mbox{for all}\ n\geq n_2.$$	
More precisely we consider the set
$$\widehat{K}_R:=(\delta_R K)\times [-L_R,L_R],$$
where $\delta_R:=\delta(R,L,K)>0$ as in Lemma \ref{lem:MD-summands-are-centered-**} (with the compact $K$ in that lemma),
$L>\Lambda_N^\prime(0)$, and $L_R>0$ as in eq. \eqref{eq:MD-summands-are-centered-ET-counting-processes}. Then we have
\begin{multline*}
P\left(\left(\frac{S_{N_n}-N_n\mu_X}{\sqrt{n/a_n}},\frac{N_n-\mathbb{E}[N_n]}{\sqrt{n/a_n}}\right)\notin\widehat{K}_R\right)\\
=P\left(\left(\left\{\frac{S_{N_n}-N_n\mu_X}{\sqrt{n/a_n}}\in\delta_RK\right\}\cap
\left\{\frac{|N_n-\mathbb{E}[N_n]|}{\sqrt{n/a_n}}\leq L_R\right\}\right)^c\right)\\
=P\left(\left\{\frac{S_{N_n}-N_n\mu_X}{\sqrt{n/a_n}}\notin\delta_RK\right\}\cup
\left\{\frac{|N_n-\mathbb{E}[N_n]|}{\sqrt{n/a_n}}>L_R\right\}\right)\\
=P\left(\frac{|N_n-\mathbb{E}[N_n]|}{\sqrt{n/a_n}}>L_R\right)+P\left(\left\{\frac{S_{N_n}-N_n\mu_X}{\sqrt{n/a_n}}\notin\delta_RK\right\}\cap
\left\{\frac{|N_n-\mathbb{E}[N_n]|}{\sqrt{n/a_n}}\leq L_R\right\}\right).
\end{multline*}
Now we estimate the summands in the final expression. We refer to eq. \eqref{eq:MD-summands-are-centered-ET-counting-processes}
to estimate the first summand. For the second one we take into account Lemma \ref{lem:MD-summands-are-centered-**} with $L=L_R>\Lambda_N^\prime(0)$
(and the definition of $\delta_R$ in this proof); then, for some $n_1\geq 1$, we have
\begin{multline*}
P\left(\left\{\frac{S_{N_n}-N_n\mu_X}{\sqrt{n/a_n}}\notin\delta_RK\right\}\cap
\left\{\frac{|N_n-\mathbb{E}[N_n]|}{\sqrt{n/a_n}}\leq L_R\right\}\right)\\
=\sum_{j=\mathbb{E}[N_n]-\sqrt{n/a_n}L_R}^{\mathbb{E}[N_n]+\sqrt{n/a_n}L_R}
P\left(\left\{\frac{S_{N_n}-N_n\mu_X}{\sqrt{n/a_n}}\notin\delta_RK\right\}\cap\{N_n=j\}\right)\\
\leq\sum_{j=0}^{nL_R}P\left(\frac{S_j-j\mu_X}{\sqrt{n/a_n}}\notin\delta_RK\right)P(N_n=j)
\leq e^{-R/a_n}\sum_{j=0}^{nL_R}P(N_n=j)\leq e^{-R/a_n};
\end{multline*}
actually we also take into account the limit in eq. \eqref{eq:limit-mean-tf-N} for the first inequality.
Thus, if we take $n_2:=\max\{n_0,n_1\}$, the statement at the beginning of the proof is proved.
This completes the proof.
\end{proof}

\subsection{MD when ``the compound sum is centered''}\label{sub:MD-centered-compound-sum}
We start with the moderate deviation result in this section. We remark that the rate function $J_{S,N}^{(2)}$ in
the next Proposition \ref{prop:MD-centered-compound-sum} can be expressed in terms of the rate function $J_{S,N}^{(1)}$ in 
Proposition \ref{prop:MD-summands-are-centered}.

\begin{proposition}\label{prop:MD-centered-compound-sum}
Consider the hypotheses in Proposition \ref{prop:MD-summands-are-centered}.
Then, for every positive sequence $\{a_n:n\geq 1\}$ such that eq. \eqref{eq:MD-conditions} holds, the sequence
$\left\{\left(\frac{S_{N_n}-\mathbb{E}[N_n]\mu_X}{\sqrt{n/a_n}},\frac{N_n-\mathbb{E}[N_n]}{\sqrt{n/a_n}}\right):n\geq 1\right\}$ satisfies
the LDP with speed $1/a_n$, and good rate function $J_{S,N}^{(2)}$ defined by
$$J_{S,N}^{(2)}(x,y):=J_{S,N}^{(1)}(x-y\mu_X,y),$$
where $J_{S,N}^{(1)}$ is the rate function in Proposition \ref{prop:MD-summands-are-centered}.
	Moreover we have the following particular cases:
	$$J_{S,N}^{(2)}(x,y):=\left\{\begin{array}{ll}
		\frac{y^2}{2\Lambda_N^{\prime\prime}(0)}&\ \mbox{if}\ x=y\mu_X\\
		\infty&\ \mbox{otherwise}
	\end{array}\right.$$
	if $\Lambda_N^\prime(0)=0$;
	$$J_{S,N}^{(2)}(x,y):=\left\{\begin{array}{ll}
		\frac{\langle x-y\mu_X,\Sigma_X^{-1}(x-y\mu_X)\rangle}{2\Lambda_N^\prime(0)}
		+\frac{y^2}{2\Lambda_N^{\prime\prime}(0)}&\ \mbox{if}\ x-y\mu_X\in\mathrm{Im}(\Sigma_X)\\
		\infty&\ \mbox{otherwise}
	\end{array}\right.$$
	if $\Lambda_N^\prime(0)>0$ and $\langle\theta,\Sigma_X\theta\rangle>0$ for all $\theta\neq\widehat{0}^*$.
\end{proposition}
\begin{proof}
The particular cases are immediate. For the remaining part of the proof we consider the function $G:B\times\mathbb{R}$ defined by
$$G(x,y):=(x+y\mu_X,y).$$
It is a continuous function and we have
$$\left(\frac{S_{N_n}-\mathbb{E}[N_n]\mu_X}{\sqrt{n/a_n}},\frac{N_n-\mathbb{E}[N_n]}{\sqrt{n/a_n}}\right)
=G\left(\left(\frac{S_{N_n}-N_n\mu_X}{\sqrt{n/a_n}},\frac{N_n-\mathbb{E}[N_n]}{\sqrt{n/a_n}}\right)\right)\quad\quad\mbox{for all}\ n\geq 1.$$
Then, by the contraction principle and Proposition \ref{prop:MD-summands-are-centered}, the desired LDP holds with good rate function 
$J_{S,N}^{(2)}$ defined by
$$J_{S,N}^{(2)}(x,y):=\inf\{J_{S,N}^{(1)}(z,w):G(z,w)=(x,y)\},$$
which coincides with the rate function in the statement of the proposition.
\end{proof}

\begin{remark}\label{rem:MD-coincidence-of-two-cases}
The rate functions in Propositions \ref{prop:MD-summands-are-centered} and \ref{prop:MD-centered-compound-sum}, i.e. $J_{S,N}^{(1)}$ and 
$J_{S,N}^{(2)}$, coincide if and only if $\mu_X=\widehat{0}$; indeed we have
$$\sum_{i=1}^{N_n}(X_i-\mu_X)=\sum_{i=1}^{N_n}X_i-\mathbb{E}[N_n]\mu_X\quad\quad\mbox{for all}\ n\geq 1$$
if and only if $\mu_X=\widehat{0}$ (when $N_n$ is not deterministic). Thus Propositions \ref{prop:MD-summands-are-centered} 
and \ref{prop:MD-centered-compound-sum} provide the same result if $\mu_X=\widehat{0}$.
\end{remark}

\begin{remark}\label{rem:MD-GET-for-centered-compound-sum}
One could try to prove Proposition \ref{prop:MD-centered-compound-sum} by following the same lines of the proof of Proposition 
\ref{prop:MD-summands-are-centered}. In such a case, if we refer again to the function $\Psi_{S,N}$ in eq. \eqref{eq:def-PsiSN}, the analogue 
of eq. \eqref{eq:MD-rf-LT-summands-are-centered} for $J_{S,N}^{(2)}$ is
\begin{equation}\label{eq:MD-rf-LT-centered-compound-sum}
J_{S,N}^{(2)}(x,y):=\sup_{\theta\in B^*,\eta\in\mathbb{R}}\{\langle\theta,x\rangle+\eta y-\Psi_{S,N}(\theta,\eta+\langle\theta,\mu_X\rangle)\}
\quad\quad\mbox{for all}\ (x,y)\in B\times\mathbb{R}.
\end{equation}
Indeed, by setting $\eta+\langle\theta,\mu_X\rangle=\xi$, we get
\begin{multline*}
J_{S,N}^{(2)}(x,y)=\sup_{\theta\in B^*,\xi\in\mathbb{R}}\{\langle\theta,x\rangle+(\xi-\langle\theta,\mu_X\rangle)y-\Psi_{S,N}(\theta,\xi)\}\\
=\sup_{\theta\in B^*,\xi\in\mathbb{R}}\{\langle\theta,x-y\mu_X\rangle+\xi y-\Psi_{S,N}(\theta,\xi)\}=J_{S,N}^{(1)}(x-y\mu_X,y),
\end{multline*}
and we can conclude as we did in the proof of Proposition \ref{prop:MD-centered-compound-sum}.
\end{remark}

\begin{remark}\label{rem:on-the-weak-convergences}
		Here we give some details on the weak convergence of 
		$\left(\frac{S_{N_n}-N_n\mu_X}{\sqrt{n}},\frac{N_n-\mathbb{E}[N_n]}{\sqrt{n}}\right)$ 
		and $\left(\frac{S_{N_n}-\mathbb{E}[N_n]\mu_X}{\sqrt{n}},\frac{N_n-\mathbb{E}[N_n]}{\sqrt{n}}\right)$ to centered Normal 
		distributions on $B\times\mathbb{R}$. Note that, in both cases, the second component is 
		$\frac{N_n-\mathbb{E}[N_n]}{\sqrt{n}}$, and it is assumed that converges weakly to a centered Normal distribution.
		In the first case we have the weak convergence to a centered Normal distribution with independent components, the covariance 
		operator for the $B$-valued component is $\Lambda_N^\prime(0)\Sigma_X$, and the variance for the real valued component is $\Lambda_N^{\prime\prime}(0)$. In the second case we have the linear transformation $G(x,y):=(x+y\mu_X,y)$ on $B\times\mathbb{R}$
		of the centered Normal distribution with independent components presented for the first case.
\end{remark}

\subsection{Some typical features}\label{sub:some-typical-features}
The aim of this section is to obtain the analogue of the limits in Remark \ref{rem:typical-features-for-N-only} (see eqs.
\eqref{eq:limit-mean-tf-N} and \eqref{eq:limit-variance-tf-N}) for the bivariate sequence in Proposition \ref{prop:LD}. 
In some sense we illustrate some typical features of moderate deviations related to the application of the G\"artner--Ellis Theorem. 
In view of what follows, we introduce some notation and preliminaries. We refer to $\mu_X$ in eq. \eqref{eq:condition-for-the-mean} and to 
$\Sigma_X$ in eq. \eqref{eq:condition-for-the-covariance}. Moreover we refer to the function $\Lambda_{S,N}$ in eq. \eqref{eq:def-LambdaSN} 
which appears in the rate function expression in Proposition \ref{prop:LD} (see $I_{S,N}$ in eq. \eqref{eq:LD-rf-LT}); we use the symbols 
$\partial_\eta$ for the standard derivative with respect to $\eta$, $D_\theta$ for the first Gateaux derivative with respect to $\theta$ 
(which is a linear functional defined on $B^*$) and $D_{\theta\theta}^{(2)}$ for the second Gateaux derivative with respect
to $\theta$ (which is a bilinear form defined on $B^*$). Then we have
$$\left.D_\theta\Lambda_X(\theta)\right|_{\theta=\widehat{0}}(v)=\langle v,\mu_X\rangle\quad\quad (\mbox{for}\ v\in B^*)$$
and
$$\left.D_{\theta\theta}^{(2)}\Lambda_X(\theta)\right|_{\theta=\widehat{0}}(u,v)=\langle u,\Sigma_Xv\rangle\quad\quad (\mbox{for}\ u,v\in B^*).$$

Now we are ready to present the limits that we want to prove.
\begin{itemize}
\item The analogue of eq. \eqref{eq:limit-mean-tf-N} consists of two limits:
\begin{equation}\label{eq:limit-mean-tf-S}
\lim_{n\to\infty}\mathbb{E}[\langle S_{N_n},v\rangle/n]=
\left.D_\theta\Lambda_{S,N}(\theta,\eta)\right|_{(\theta,\eta)=(\widehat{0},0)}(v)\quad\quad (\mbox{for}\ v\in B^*),
\end{equation}
and eq. \eqref{eq:limit-mean-tf-N} itself.
\item The analogue of eq. \eqref{eq:limit-variance-tf-N} consists of three limits:
\begin{equation}\label{eq:limit-cov-tf-S,S}
\lim_{n\to\infty}n\mathrm{Cov}\left(\frac{\langle S_{N_n},u\rangle}{n},\frac{\langle S_{N_n},v\rangle}{n}\right)=
\left.D_{\theta\theta}^{(2)}\Lambda_{S,N}(\theta,\eta)\right|_{(\theta,\eta)=(\widehat{0},0)}(u,v)\quad\quad (\mbox{for}\ u,v\in B^*),
\end{equation}
\begin{equation}\label{eq:limit-cov-tf-N,S}
\lim_{n\to\infty}n\mathrm{Cov}\left(\frac{N_n}{n},\frac{\langle S_{N_n},v\rangle}{n}\right)=
\left.\partial_\eta D_\theta\Lambda_{S,N}(\theta,\eta)\right|_{(\theta,\eta)=(\widehat{0},0)}(v)\quad\quad (\mbox{for}\ v\in B^*),
\end{equation}
and eq. \eqref{eq:limit-variance-tf-N} itself.
\end{itemize}

\begin{remark}\label{rem:typical-features}
If we consider $\Psi_{S,N}(\theta,\eta+\langle\theta,\mu_X\rangle)$ which appears in eq. \eqref{eq:MD-rf-LT-centered-compound-sum}, after
some computations we can say that
$$\Psi_{S,N}(\theta,\eta+\langle\theta,\mu_X\rangle)=\frac{1}{2}\left(b_1(\theta)+2b_2(\theta)\eta+\Lambda_N^{\prime\prime}(0)\eta^2\right)
\quad\quad\mbox{for all}\ (\theta,\eta)\in B^*\times\mathbb{R},$$
where
$$b_1(\theta):=\Lambda_N^{\prime\prime}(0)(\langle\theta,\mu_X\rangle)^2+\Lambda_N^\prime(0)\langle\theta,\Sigma_X\theta\rangle\quad\quad\mbox{and}
\quad\quad b_2(\theta):=\Lambda_N^{\prime\prime}(0)\langle\theta,\mu_X\rangle.$$
Then $b_1(\theta)$ coincides with the limit value in eq. \eqref{eq:limit-cov-tf-S,S} with $u=v=\theta$ (see also eq. \eqref{eq:limit-cov-tf-S,S-explicit}) 
and $b_2(\theta)$ coincides with the limit value in eq. \eqref{eq:limit-cov-tf-N,S} with $v=\theta$ (see also eq. \eqref{eq:limit-cov-tf-N,S-explicit});
moreover we recall the limit in eq. \eqref{eq:limit-variance-tf-N} for the coefficient $\Lambda_N^{\prime\prime}(0)$ of $\eta^2$.
\end{remark}

We start with more explicit expressions of the limit values in eqs. \eqref{eq:limit-mean-tf-S}, \eqref{eq:limit-cov-tf-S,S} and 
\eqref{eq:limit-cov-tf-N,S}. Firstly we take $v\in B^*$ and we have
$$D_\theta\Lambda_{S,N}(\theta,\eta)(v)=\Lambda_N^\prime(\eta+\Lambda_X(\theta))D_\theta\Lambda_X(\theta)(v)$$
and
$$\partial_\eta D_\theta\Lambda_{S,N}(\theta,\eta)(v)=\Lambda_N^{\prime\prime}(\eta+\Lambda_X(\theta))D_\theta\Lambda_X(\theta)(v);$$
therefore we get
$$\left.D_\theta\Lambda_{S,N}(\theta,\eta)\right|_{(\theta,\eta)=(\widehat{0},0)}(v)=\Lambda_N^\prime(0)\langle v,\mu_X\rangle$$
for the limit value in eq. \eqref{eq:limit-mean-tf-S}, and
\begin{equation}\label{eq:limit-cov-tf-N,S-explicit}
\left.\partial_\eta D_\theta\Lambda_{S,N}(\theta,\eta)\right|_{(\theta,\eta)=(\widehat{0},0)}(v)=\Lambda_N^{\prime\prime}(0)\langle v,\mu_X\rangle
\end{equation}
for the limit value in eq. \eqref{eq:limit-cov-tf-N,S}. Moreover we take $u,v\in B^*$ and we have
$$D_{\theta\theta}^{(2)}\Lambda_{S,N}(\theta,\eta)(u,v)
=\Lambda_N^{\prime\prime}(\eta+\Lambda_X(\theta))D_\theta\Lambda_X(\theta)(u)D_\theta\Lambda_X(\theta)(v)
+\Lambda_N^\prime(\eta+\Lambda_X(\theta))D_{\theta\theta}^{(2)}\Lambda_X(\theta)(u,v);$$
therefore we get
\begin{equation}\label{eq:limit-cov-tf-S,S-explicit}
\left.D_{\theta\theta}^{(2)}\Lambda_{S,N}(\theta,\eta)\right|_{(\theta,\eta)=(\widehat{0},0)}(u,v)=
\Lambda_N^{\prime\prime}(0)\langle u,\mu_X\rangle\langle v,\mu_X\rangle+\Lambda_N^\prime(0)\langle u,\Sigma_Xv\rangle
\end{equation}
for the limit value in eq. \eqref{eq:limit-cov-tf-S,S}.

Now some details on the computations of the limits. The limit in eq. \eqref{eq:limit-mean-tf-S} can be easily obtained noting that
$$\mathbb{E}[\langle S_{N_n},v\rangle/n]=\mathbb{E}[N_n/n]\mathbb{E}[\langle X_1,v\rangle]
=\mathbb{E}[N_n/n]\langle v,\mu_X\rangle,$$
and by taking into account the limit in eq. \eqref{eq:limit-mean-tf-N}.
The limit in eq. \eqref{eq:limit-cov-tf-N,S} can be easily obtained noting that (after some easy computations)
$$n\mathrm{Cov}\left(\frac{N_n}{n},\frac{\langle S_{N_n},v\rangle}{n}\right)=n\mathrm{Var}[N_n/n]\langle v,\mu_X\rangle,$$
and by taking into account the limit in eq. \eqref{eq:limit-variance-tf-N}.
The limit in eq. \eqref{eq:limit-cov-tf-S,S} can be obtained noting that
\begin{multline*}
n\mathrm{Cov}\left(\frac{\langle S_{N_n},u\rangle}{n},\frac{\langle S_{N_n},v\rangle}{n}\right)
=n\left(\mathbb{E}[N_n/n^2]\mathbb{E}[\langle X_1,u\rangle\langle X_1,v\rangle]\right.\\
\left.+\mathbb{E}[N_n(N_n-1)/n^2]\mathbb{E}[\langle X_1,u\rangle]\mathbb{E}[\langle X_1,v\rangle]
-\mathbb{E}^2[N_n/n]\mathbb{E}[\langle X_1,u\rangle]\mathbb{E}[\langle X_1,v\rangle]\right)\\
=\mathbb{E}[N_n/n]\mathrm{Cov}(\langle X_1,u\rangle,\langle X_1,v\rangle)
+n\mathrm{Var}[N_n/n]\mathbb{E}[\langle X_1,u\rangle]\mathbb{E}[\langle X_1,v\rangle]\\
=\mathbb{E}[N_n/n]\langle u,\Sigma_Xv\rangle
+n\mathrm{Var}[N_n/n]\langle u,\mu_X\rangle\langle v,\mu_X\rangle,
\end{multline*}
and by taking into account the limits in eqs. \eqref{eq:limit-mean-tf-N} and \eqref{eq:limit-variance-tf-N}.

\paragraph{A particular case: $B=C(K,\mathbb{R})$, where $K$ is a compact set.}
Here we consider the case in which $B$ is the separable Banach space formed by the real-valued continuous function 
defined on a compact set $K$. In this case $B^*$ is the family of signed measures on $K$ and, for every $f\in B$ and 
$v\in B^*$, we have
$$\langle f,v\rangle=\int_K f(s)dv(s).$$
Then eq. \eqref{eq:condition-for-the-mean} holds with $\mu_X(\cdot)=\mathbb{E}[X_1(\cdot)]\in B$ because we have
$$\mathbb{E}[\langle X_1,v\rangle]=\int_K\mathbb{E}[X_1(s)]dv(s)$$
for every $v\in B^*$. The statement in eq. \eqref{eq:condition-for-the-covariance} holds because, for every $v\in B^*$,
we can consider $\Sigma_Xv\in B=C(K,\mathbb{R})$ defined by
$$\Sigma_Xv(s):=\int_K\mathrm{Cov}(X_1(s),X_1(t))dv(t)\quad\quad (s\in K),$$
and we have
$$\langle u,\Sigma_Xv\rangle=\langle\Sigma_Xu,v\rangle=\int_{K\times K}\mathrm{Cov}(X_1(s),X_1(t))du(s)dv(t)
=\mathrm{Cov}(\langle X_1,u\rangle,\langle X_1,v\rangle)$$
for every $u,v\in B^*$.

\section{Examples}\label{sec:examples}
In this section we present some examples. We start with some examples for the sequence $\{N_n:n\geq 1\}$; in particular, we refer to 
Conditions \ref{cond:LD-*} and \ref{cond:MD-*}, and to the limits in Remark \ref{rem:typical-features-for-N-only} (see eqs. 
\eqref{eq:limit-mean-tf-N} and \eqref{eq:limit-variance-tf-N}). Finally we present an example for the sequence $\{X_n:n\geq 1\}$.

\subsection{Some examples for the sequence $\{N_n:n\geq 1\}$}\label{sub:ex-N}
We start with Examples \ref{ex:random-walk} and \ref{ex:NHPP} and we omit the details; we only say that we can refer to the proof of 
Theorem 3.7.1 in \cite{DZ:97} (case $d=1$) for the limit in eq. \eqref{eq:MD-limit-for-marginal-N} for Example \ref{ex:random-walk}.

\begin{example}\label{ex:random-walk}
	Let $\{Z_n:n\geq 1\}$ be a sequence of i.i.d. nonegative integer valued random variables with moment generating function 
	$\mathcal{G}_Z(\eta):=\mathbb{E}[e^{\eta Z_1}]$ assumed to be finite for every $\eta\in\mathbb{R}$.
	Then, if we set
	$$N_n:=Z_1+\cdots+Z_n\quad\quad\mbox{for every}\ n\geq 1,$$
    we have $\Lambda_N(\eta)=\log\mathcal{G}_Z(\eta)$ for every $\eta\in\mathbb{R}$, $\Lambda_N^\prime(0)=\mathbb{E}[Z_1]$ and 
    $\Lambda_N^{\prime\prime}(0)=\mathrm{Var}[Z_1]$.
\end{example}

\begin{example}\label{ex:NHPP}
	Let $\{M(t):t\geq 0\}$ be a homogeneous Poisson process, and let $t\mapsto\lambda(t)$ be a locally integrable intensity
	(i.e. a a locally integrable positive function) such that
	$$\lim_{t\to\infty}\frac{\int_0^t\lambda(s)ds}{t}=\widetilde{\lambda}\in(0,\infty).$$
	Then, if we set
	$$N_n:=M\left(\int_0^n\lambda(s)ds\right)\quad\quad\mbox{for every}\ n\geq 1,$$
	we have $\Lambda_N(\eta)=\widetilde{\lambda}(e^\eta-1)$ for every $\eta\in\mathbb{R}$, $\Lambda_N^\prime(0)=\widetilde{\lambda}$ and
	$\Lambda_N^{\prime\prime}(0)=\widetilde{\lambda}$.
\end{example}

We can also present other examples. The first one concerns the so-called alternative fractional Poisson process; see the family of the 
random variables presented in \cite{BeghinMacciSPL2013} as a suitable deterministic time-change of the one defined by eq. (4.1) in 
\cite{BeghinOrsingherEJP2009}.

\begin{example}\label{ex:FPP-alternative}
	Let $\nu\in(0,1]$ and $\lambda>0$ be arbitrarily fixed. Then, for every $n\geq 1$, by a suitable modified version of eq. (4.4) in 
	\cite{BeghinOrsingherEJP2009} we have
	$$\mathbb{E}[e^{\eta N_n}]=\frac{E_{\nu,1}(e^\eta\lambda n^\nu)}{E_{\nu,1}(\lambda n^\nu)}\quad\quad\mbox{for every}\ \eta\in\mathbb{R},$$
	where $E_{\nu,\beta}(x):=\sum_{r\geq 0}\frac{x^r}{\Gamma(\nu r+\beta)}$ is the two-parametric Mittag-Leffler function (see, e.g., 
	\cite{GorenfloKilbasMainardiRogosin2014}). Then we have $\Lambda_N(\eta)=\lambda^{1/\nu}(e^{\eta/\nu}-1)$ for every $\eta\in\mathbb{R}$, 
	$\Lambda_N^\prime(0)=\frac{\lambda^{1/\nu}}{\nu}$ and $\Lambda_N^{\prime\prime}(0)=\frac{\lambda^{1/\nu}}{\nu^2}$. Here we give a list
	of references of interest: the proof of Proposition 4.1 (and in particular the function $\Lambda_N$ in Condition 
	\ref{cond:LD-*}) in \cite{BeghinMacciSPL2013}; the proof of Proposition 2 (for the case $m=1$) in \cite{BeghinMacciSPL2017}
	for the statements in Condition \ref{cond:MD-*}; the equality
	$$\mathbb{E}[N_n]=\frac{\lambda n^\nu}{\nu}\frac{E_{\nu,\nu}(\lambda n^\nu)}{E_{\nu,1}(\lambda n^\nu)}$$
	(which is a suitable modified version of eq. (4.6) in \cite{BeghinOrsingherEJP2009}) for the limit in eq. \eqref{eq:limit-mean-tf-N}, 
	indeed this limit can be checked by taking into account the asymptotic behavior of $E_{\nu,\beta}(x)$ as $x$ goes to infinity (see, e.g.,
	eq. (4.4.16) in \cite{GorenfloKilbasMainardiRogosin2014}).
\end{example}

A further example concerns sums of independent Bernoulli random variables. In this way we consider a generalization of the results presented in
\cite{GiulianoMacciPacchiarottiJSPI2015} (Sections 3.1 and 3.2) for runs. We recover that case by taking $\mathcal{X}=[0,1]$, 
$p(x):=e^{-\lambda cx}$ for $\lambda,c>0$, and $x_{j,n}:=(j-1)\varepsilon_n$ for $1\leq j\leq n$, where $\varepsilon_n\sim\frac{1}{n}$; therefore
the weak convergence in eq. \eqref{eq:weak-convergence-in-runs-generalization} holds (see below) and $R$ is the Lebesgue measure on $\mathcal{X}$.

\begin{example}\label{ex:runs-generalization}
	Let $\mathcal{X}$ be a topological space, let $R$ be a positive measure on $\mathcal{X}$, and let $p:\mathcal{X}\to[0,1]$ be a continuous
	function. Moreover let $\{x_{j,n}:1\leq j\leq n\}$ be a family of points in $\mathcal{X}$ such that
	\begin{equation}\label{eq:weak-convergence-in-runs-generalization}
	    \frac{1}{n}\sum_{j=1}^n1_{\{x_{j,n}\in\cdot\}}\ \mbox{converges weakly to}\ R(\cdot)\quad\quad\mbox{as}\ n\to\infty,
	\end{equation}
	and let $\{Z_{j,n}:1\leq j\leq n\}$ be a sequence of random variables such that, for every $n\geq 1$, $Z_{1,n},\ldots,Z_{n,n}$ are 
	independent and
	$$P(Z_{j,n}=1)=1-P(Z_{j,n}=0)=p(x_{j,n})\quad\quad\mbox{for}\ 1\leq j\leq n.$$
    Then, if we set
    $$N_n:=Z_{1,n}+\cdots+Z_{n,n}\quad\quad\mbox{for every}\ n\geq 1,$$
    we have $\Lambda_N(\eta)=\int_{\mathcal{X}}\log(1+p(x)(e^\eta-1))dR(x)$ for every $\eta\in\mathbb{R}$, 
    $\Lambda_N^\prime(0)=\int_{\mathcal{X}}p(x)dR(x)$ and $\Lambda_N^{\prime\prime}(0)=\int_{\mathcal{X}}p(x)(1-p(x))dR(x)$. It is
    possible to check all these statements by adapting all the arguments and computations in \cite{GiulianoMacciPacchiarottiJSPI2015}.
\end{example}

A final example concerns renewal processes.

\begin{example}\label{ex:renewal-process}
	Let $\{T_n:n\geq 1\}$ be a sequence of i.i.d. positive random variables. We set $\kappa_T(\eta):=\log\mathbb{E}[e^{\eta T_1}]$ 
	and we assume that it is finite in a neighbourhood of the origin $\eta=0$, and $\kappa_T(-\infty)=-\infty$. In
	order to avoid trivialities we also assume that $(\kappa_T)^\prime(0)$ and $(\kappa_T)^{\prime\prime}(0)$ are finite and
	positive. Then, if we set
	$$N_n:=\sum_{k\geq 1}1_{\{T_1+\cdots+T_k\leq n\}},\quad\quad\mbox{for every}\ n\geq 1,$$
	the function $\Lambda_N$ in Condition \ref{cond:LD-*} is well-defined and we have
	$\Lambda_N(\eta):=-(\kappa_T)^{-1}(-\eta)$ for every $\eta\in\mathbb{R}$. Note that the function $\Lambda_N$ is 
    finite everywhere because $\kappa_T(-\infty)=-\infty$ (see eqs. (12) and (13) in Theorem 1 in \cite{GlynnWhittQS1994});
    moreover 
    $\Lambda_N^\prime(0)=1/(\kappa_T)^\prime(0)$ and $\Lambda_N^{\prime\prime}(0)=(\kappa_T)^{\prime\prime}(0)/((\kappa_T)^\prime(0))^3$.
    Then, if one can prove Condition \ref{cond:MD-*} and the limit in eq. \eqref{eq:limit-mean-tf-N}, we can apply the results in 
    this paper. Note that this happens if the random variables $\{T_n:n\geq 1\}$ are exponentially distributed, and therefore 
    $\{N_n:n\geq 1\}$ is the discrete time version of the Poisson process; in particular, for $\lambda>0$, we have
    $$\kappa_T(\eta)=\log\frac{\lambda}{\lambda-\eta}\quad\quad\mbox{for}\ \eta<\lambda$$
    and $\kappa_T(\eta)=\infty$ otherwise, and $\Lambda_N(\eta)=\lambda(e^\eta-1)$ for all $\eta\in\mathbb{R}$.
    For completeness we recall that $\frac{N_n-\mathbb{E}[N_n]}{\sqrt{n}}$ converges weakly to a centered Normal 
    distribution with variance $\Lambda_N^{\prime\prime}(0)$ (see, e.g., Problem 27.15 in \cite{Billingsley-Probability-and-Measure-1995}).
\end{example}

\subsection{An example for the sequence $\{X_n:n\geq 1\}$}\label{sub:ex-S}
One could present several examples for $\{X_n:n\geq 1\}$ with well-known formulas. For instance, if $B=C(K,\mathbb{R}^h)$ for a compact
set $K$, we could consider i.i.d. continuous Gaussian processes or diffusions $\{X_n:n\geq 1\}$ on some $\mathbb{R}^h$. Another example 
consists of i.i.d. random variables $\{X_n:n\geq 1\}$ with Normal distribution on a Banach space $B$. In our opinion the following 
example has some more interesting aspects. In particular we discuss some formulas in Propositions \ref{prop:LD} and
\ref{prop:MD-summands-are-centered}.

\begin{example}\label{ex:finite-valued-X}
	Let $B$ be a Banach space such that $\frac{S_n-n\mu_X}{\sqrt{n}}$ converges weakly to a centered normal distribution,
	as $n\to\infty$ (see Condition \ref{cond:MD-*}); the interested reader can refer to \cite{CLT_1976} and the references cited therein.
	Moreover, let $\{X_n:n\geq 1\}$ be such that
	$$P(X_1=u_i)=p_i\quad\quad\mbox{for}\ 1\leq i\leq m$$
	for some $m\geq 2$, where $p_1,\ldots,p_m>0$ and $p_1+\cdots+p_m=1$, and $u_1,\ldots,u_m\in B$ are linearly independent.
\end{example}

We start with the formulas of Proposition \ref{prop:LD}. We have
$$\Lambda_X(\theta)=\log\left(\sum_{i=1}^me^{\langle\theta,u_i\rangle}p_i\right)\quad\quad \mbox{for all}\ \theta\in B^*.$$
It is useful to consider the set $\mathcal{S}$ of linear convex combinations of $u_1,\ldots,u_m$; thus $x\in\mathcal{S}$
if and only if
\begin{equation}\label{eq:linear-convex-combination}
	x=c_1u_1+\cdots c_mu_m,\quad\quad \mbox{where}\ c_1,\ldots,c_m\geq 0\ \mbox{and}\ c_1+\cdots+c_m=1.
\end{equation}
Then one can check that
$$\Lambda_X^*(x)=\left\{\begin{array}{ll}
	\sum_{i=1}^mc_i\log\left(\frac{c_i}{p_i}\right)&\ \mbox{if}\ x\in\mathcal{S}\ (\mbox{and}\ \eqref{eq:linear-convex-combination}\ \mbox{holds})\\
	\infty&\ \mbox{otherwise},
\end{array}\right.$$
where $0\log 0=0$. This can be checked as follows.
\begin{itemize}
	\item If eq. \eqref{eq:linear-convex-combination} holds with $c_1,\ldots,c_m>0$, then
	$$\sup_{\theta\in B^*}\left\{\langle\theta,x\rangle-\log\left(\sum_{i=1}^me^{\langle\theta,u_i\rangle}p_i\right)\right\}$$
	is attained at any arbitrary solution $\theta\in B^*$ of the equation $x=\frac{\sum_{i=1}^me^{\langle\theta,u_i\rangle}p_iu_i}{\sum_{i=1}^me^{\langle\theta,u_i\rangle}p_i}$;
	then we have $\langle\theta,u_i\rangle=\log\left(\frac{c_i}{p_i}\right)$ for all $1\leq i\leq m$, and we easily get the
	desired expression of $\Lambda_X^*(x)$.
	\item If eq. \eqref{eq:linear-convex-combination} holds with $c_i=0$ for some $1\leq i\leq m$ (possibly not unique), we easily
	get the desired expression by the continuity of $\Lambda_X^*(x)$ when $x$ belongs to the interior of $\mathcal{S}$.
	\item If $x\notin\mathcal{S}$ we get $\Lambda_X^*(x)=\infty$ by taking into account the lower bound for open sets provided 
	by the Cramér Theorem in Banach spaces (see, e.g., \cite{deAcosta2}).
\end{itemize}

Now we consider the formulas of Proposition \ref{prop:MD-summands-are-centered}. In order to avoid trivialities we assume that 
$\Lambda_N^\prime(0)\neq 0$. It is useful to consider the set $\mathcal{L}$ of linear combinations of $u_1,\ldots,u_m$ such that the 
sum of the coefficients are equal to zero; thus $x\in\mathcal{L}$ if and only if
\begin{equation}\label{eq:linear-combination-sum-zero}
	x=c_1u_1+\cdots c_mu_m,\quad\quad \mbox{where}\ c_1,\ldots,c_m\in\mathbb{R}\ \mbox{and}\ c_1+\cdots+c_m=0.
\end{equation}
Then
$$\sup_{\theta\in B^*}\left\{\langle\theta,x\rangle-\frac{\Lambda_N^\prime(0)}{2}\langle\theta,\Sigma_X\theta\rangle\right\}
=\left\{\begin{array}{ll}
	\frac{1}{2\Lambda_N^\prime(0)}\sum_{j=1}^{m-1}c_j\left(\frac{c_j}{p_j}-\frac{c_m}{p_m}\right)&\ \mbox{if}\ x\in\mathcal{L}
	\ (\mbox{and}\ \eqref{eq:linear-combination-sum-zero}\ \mbox{holds})\\
	\infty&\ \mbox{otherwise}.
\end{array}\right.$$
This can be checked as follows.
\begin{itemize}
	\item If eq. \eqref{eq:linear-combination-sum-zero} holds, then
	$$\sup_{\theta\in B^*}\left\{\langle\theta,x\rangle-\frac{\Lambda_N^\prime(0)}{2}\langle\theta,\Sigma_X\theta\rangle\right\}$$
	is attained at any arbitrary solution $\theta\in B^*$ of $x=\Lambda_N^\prime(0)\Sigma_X\theta$; then we have
	$$\left(\langle\theta,u_j\rangle-\sum_{h=1}^m\langle\theta,u_h\rangle p_h\right)p_j=\frac{c_j}{\Lambda_N^\prime(0)}
	\quad\quad \mbox{for}\ j\in\{1,\ldots,m-1\},$$
	which yields (after some computations)
	$$\langle\theta,u_j\rangle=\frac{1}{\Lambda_N^\prime(0)}\left(\frac{c_j}{p_j}-\frac{c_m}{p_m}\right)+\langle\theta,u_m\rangle
	\quad\quad \mbox{for}\ j\in\{1,\ldots,m-1\}.$$
	In conclusion we easily get the desired expression of $\sup_{\theta\in B^*}\left\{\langle\theta,x\rangle-\frac{\Lambda_N^\prime(0)}{2}\langle\theta,\Sigma_X\theta\rangle\right\}$ noting that
	$$\langle\theta,x\rangle-\frac{\Lambda_N^\prime(0)}{2}\langle\theta,\Sigma_X\theta\rangle=
	\frac{\langle\theta,x\rangle}{2}=\frac{1}{2}\sum_{j=1}^mc_j\langle\theta,u_j\rangle.$$
	\item If $x\notin\mathcal{L}$ we can adapt the argument above for $x\notin\mathcal{S}$ to say that $\Lambda_X^*(x)=\infty$.
\end{itemize}

\paragraph{Acknowledgements.}
The idea to study the problems addressed in this article was inspired by some discussions of one of the authors with Nikolai Leonenko.
The authors thank two referees for some useful comments, and Carlo Sinestrari for some useful hints on the content of reference 
\cite{Rockafeller}.

\paragraph{Funding.}
The authors are supported by MIUR Excellence Department Project awarded to the Department of Mathematics, University
of Rome Tor Vergata (CUP E83C23000330006), by University of Rome Tor Vergata (project ``Asymptotic Properties in Probability''
(CUP E83C22001780005)) and by Indam-GNAMPA.

\bibliography{bibbase}
\bibliographystyle{plain}

\end{document}